\documentclass{article}

\usepackage{mathrsfs} 
\usepackage{amsthm, amsmath, amsfonts, amssymb}
\usepackage{dsfont}
\usepackage{hyperref}
\hypersetup{
  linkcolor=green,
  citecolor=green,
  colorlinks=false,
}

\newtheorem{define}{Definition}[section]
\newtheorem{thm}[define]{Theorem}
\newtheorem{lemma}[define]{Lemma}

\newcommand{\func}{f}
\newcommand{\funcc}{g}
\newcommand{\funccc}{h}
\newcommand{\funcs}{\mathscr{F}}

\begin{document}

\title{Blocker size via matching minors}
\author{Nikola Yolov}
\maketitle
\begin{abstract}
  Finding the maximum number of maximal independent sets
  in an $n$-vertex graph $G$, $i(G)$, from a restricted class
  is an extensively studied problem.
  Let $kK_2$ denote the matching of size $k$,
  that is a graph with $2k$ vertices and $k$ disjoint edges.
  A graph with an induced copy of $kK_2$
  contains at least $2^k$ maximal independent sets.
  The other direction was established in a series of papers
  \cite{farber1, balas_yu, farber2}
  finally yielding $i(G) \le (n/k)^{2k}$ for a graph $G$
  without an induced $(k+1)K_2$.
  Alekseev proved in \cite{alekseev}
  that $i(G)$ is at most the number of induced matchings of $G$.

  This work generalises the aforementioned results to clutters.
  The right substructures in this setting are minors rather than
  induced subgraphs.
  Maximal independent sets of a clutter $\mathcal{H}$
  are in one-to-one correspondence to the sets of its blocker, $b(\mathcal{H})$,
  hence $i(\mathcal{H}) = |b(\mathcal{H})|$.
  We show that
  \[
  |b(\mathcal{H})| \le
  \sum_{m=0}^{k \cdot f(r)}{|\mathcal{H}| \choose m} {r \choose 2}^m
  \]
  for a $(k+1)K_2$-minor-free clutter $\mathcal{H}$
  where $f(r) = (2r-3)2^{r-2}$ and $r$ is the maximum size of a set in $\mathcal{H}$.
  A key step in the proofs is,  similarly to Alekseev's result,
  showing that $i(\mathcal{H})$ is bounded
  by the number of a substructure called semi-matching,
  and then proving a dependence between the number of semi-matchings
  and the number of minor matchings.
  Note that similarly to graphs, a clutter containing a $kK_2$ minor
  has at least $2^k$ maximal independent sets.

  From a computational perspective,
  a polynomial number of independent sets is particularly interesting.
  Our results lead to polynomial algorithms for restricted instances
  of many problems including Set Cover and k-SAT.
\end{abstract}

\section{Introduction}
\subsection{Previous work}
An \emph{independent (or stable) set} of a graph $G$ is a subset of its vertices
not containing an edge.
Let $i(G)$ be the number of inclusion maximal independent sets of $G$.
The first result on $i(G)$ is by Moon and Moser \cite{moon1965cliques}
showing that $i(G) \le 3^{n/3}$ for an $n$-vertex graph $G$.
The exact number has been found for
general graphs,
connected graphs,
trees,
forests,
triangle-free graphs
and other classes.
This research sparked the interested of finding
a family of classes controlled by a parameter
where the maximum number of maximal independent sets
can be bounded in terms of the parameter.
It is easy to see that a graph $G$ with an induced matching of size $k$
contains at least $2^k$ maximal independent sets.
Farber proved in \cite{farber1} that an $n$-vertex graph $G$ without an
induced $2K_2$ has $O(n^2)$ maximal independent sets.
Let $\rho = \rho(G) = |\{uv \in E(G) : N(u) \cup N(v) \neq V(G)\}|$,
and let $\nu = \nu(G)$ be the maximal $k$ such that $G$ contains $kK_2$
as an induced graph.
Balas and Yu generalised the result of Farber in \cite{balas_yu} showing that
$i(G) \le \rho^\nu + 1$.
Farber, Hujter and Tuza proved in \cite{farber2}
that $i(G) < \sqrt[3] \nu e^{1-\nu} {n \choose 2\nu}$
for a graph on $n \ge 4 \nu$ vertices.
All these results
are generalised by an elegant theorem by Alekseev \cite{alekseev}
stating that $i(G)$ is at most the number of induced matchings of $G$.

A number of otherwise NP-complete problems can be solved in polynomial time if
the list of all maximal stable sets of the underlying graph is given as input,
as a trivial example consider the maximum independent set problem.
Other problems have their complexity reduced.
For instance, given this input, the chromatic number problem can be
reduced to the set cover problem, which is $\log n$-approximable,
while the chromatic number is not $n^{1-\epsilon}$ approximable
for any $\epsilon > 0$, assuming  NP$\not\subseteq$ZPP\cite{hardness_chromatic}.
The set of maximal independent sets of a graph can be computed
in time polynomial with respect to their number \cite{gen_stable},
and hence graphs with polynomial-sized $i(G)$ are of a particular interest
for computer science.
A class of graphs $\mathcal{G}$ has polynomially bounded $i(G)$
if there is a polynomial $p(x)$ such that $i(G) \le p(v(G))$
for every $G \in \mathcal{G}$.
The results above imply that a hereditary class of graphs $\mathcal{G}$,
that is a class closed under taking induced subgraphs,
has polynomially bounded $i(G)$ if and only if $kK_2 \notin \mathcal{G}$
for some $k$.

\subsection{New results}
The aim of this paper is to generalise these results to hypergraphs.
An \emph{independent set} of a hypergraph $\mathcal{H}$
is a subset of $V(\mathcal{H})$ containing no $S \in \mathcal{H}$.
Note that if a hypergraph $H$ contains two hyperedges $g$ and $h$,
such that $g \subset h$,
in which case we say that $h$ is \emph{subsumed} by $g$,
the set of maximal independent sets of $H$ does not change
if $h$ is removed.
Therefore subsumed edges play no role here,
and hence we lose nothing by assuming there are no such edges.
This leads us to the first definition.

\begin{define} [Clutter]
  A \emph{set system}
  is a finite set of finite subsets of the natural numbers $\mathbb{N}$.
  A set system is called \emph{clutter}, \emph{Sperner family}
  or \emph{antichain} if its sets are incomparable.
  Every set system $\mathcal{H}$ we associate with a vertex set
  $V(\mathcal{H}) := \bigcup_{S \in \mathcal{H}} S$
  and rank $rk(\mathcal{H}) := \max_{S \in \mathcal{H}}|S|$.
  Given a set system $\mathcal{H}$ we define $cl(\mathcal{H})$ to be the clutter
  obtained from $\mathcal{H}$ by removing all non-minimal (subsumed) edges.
\end{define}
\noindent
The \emph{edge clutter} of a graph $G$
is the clutter composed of the edges of $G$.

In the context of clutters, the complements of independent sets,
called \emph{transversals}, are more natural to work with.

\begin{define} [Blocker]
  A \emph{transversal} of a clutter $\mathcal{H}$
  is a set intersecting every $S \in \mathcal{H}$.
  A \emph{minimal transversal} is a transversal not containing
  another transversal as a subset.
  The \emph{blocker} of $\mathcal{H}$, denoted $b(\mathcal{H})$,
  is the clutter consisting of all minimal transversals of $\mathcal{H}$.
\end{define}

One reason why transversals are preferred
is because $b(b(\mathcal{H})) = \mathcal{H}$ for every clutter $\mathcal{H}$,
and this property is not shared with independent sets.
For example, take $C_6 = \{12, 23, 34, 45, 56, 61\}$
and note that
$indep(C_6) = \{14, 25, 36, 135, 246\}$,
$indep(indep(C_6)) = \{123, 234, 345, 456, 561, 612\} \neq C_6$.

\begin{define} [Deletion and Contraction]
  Suppose $\mathcal{H}$ is a clutter and $v \in \mathbb{N}$.
  Define
  \begin{align*}
    \mathcal{H} \backslash v &:= \{S \in  \mathcal{H}: v \notin S\}
    \text{ to be $\mathcal{H}$ with $v$ deleted, and}\\
    \mathcal{H} / v &:= cl(\{S - v : S \in \mathcal{H}\})
    \text{ to be $\mathcal{H}$ with $v$ contracted.}
  \end{align*}
\end{define}

\begin{define} [Minor]
  We say that $\mathcal{F}$ is a \emph{minor} of $\mathcal{H}$
  if $\mathcal{F}$ can be obtained from $\mathcal{H}$
  through a series of deletions and contractions.
  We write $\mathcal{F} \subseteq_m \mathcal{H}$
  to denote that $\mathcal{F}$ is isomorphic to a minor of $\mathcal{H}$
  and say that $\mathcal{H}$ is $\mathcal{F}$-minor-free if
  $\mathcal{F} \not\subseteq_m \mathcal{H}$.
\end{define}
\noindent
Note that $\mathcal{H} \backslash v = cl(\mathcal{H} \backslash v)$,
and hence all minors of a clutter are clutters.

We are in position to state our main theorem.
\begin{thm}
  \label{thm:main}
  Suppose $\mathcal{H}$ is a $(k+1)K_2$-minor-free rank $r$ clutter.
  Then
  \[
  |b(\mathcal{H})| \le \sum_{m=0}^{k \cdot (2r-3)2^{r-2}}
  {|\mathcal{H}| \choose m} {r \choose 2}^m.
  \]
\end{thm}
\noindent
The theorem holds with equality for $\mathcal{H} \cong kK_2$.
Similarly to the case of graphs,
we see that a minor-closed class of clutters has polynomially bounded
$|b(\mathcal{H})|$ if and only if it does not contain $kK_2$ for some $k$.

\section{Notation and preliminaries}
Denote the set $\{1, \ldots n\}$ by $[n]$.
It is often convenient to write $S \cup x$ and $S - x$
instead of $S \cup \{x\}$ and $S \setminus \{x\}$.

We quickly revise some well-known properties of clutters and blockers
used throughout the proofs.
\begin{define} [Join and Meet]
  Suppose $\mathcal{H}$ and $\mathcal{F}$ are clutters.
  We define
  \begin{align*}
    \mathcal{H} \vee \mathcal{F} &:= cl(\mathcal{H} \cup \mathcal{F})
    \text{ to be the \emph{join} of } \mathcal{H} \text{ and } \mathcal{F}
    \text{, and}\\
    \mathcal{H} \wedge \mathcal{F} &:=
    cl(\{S \cup T : S \in \mathcal{H}, T \in \mathcal{F}\})
    \text{ to be the \emph{meet} of } \mathcal{H} \text{ and } \mathcal{F}.
  \end{align*}
\end{define}

\begin{lemma} [Algebraic properties]
  \label{thm:algebra}
  Let $\mathcal{H}$, $\mathcal{G}$ and $\mathcal{F}$ be clutters
  and let $v$ and $u$ be distinct elements of $\mathbb{N}$.
  The operations defined above have the following properties:
  \begin{enumerate}
  \item
    deletion and contraction commute, that is
    $\mathcal{H} \backslash v \backslash u = \mathcal{H} \backslash u \backslash v$,
    $\mathcal{H} \backslash v / u = \mathcal{H} / u\backslash v$ and
    $\mathcal{H} / v / u = \mathcal{H} / u / v$;
  \item
    if we denote the set of all clutters by $CL(\mathbb{N})$,
    the clutters $\emptyset$ and $\{\emptyset\}$
    by $\widehat{0}$ and $\widehat{1}$ respectively,
    then $(CL(\mathbb{N}), \vee, \wedge, \widehat{0}, \widehat{1})$
    forms a bounded distributive lattice, that is
    \begin{center}
      \begin{tabular}{ c c c}
        Commutative laws &
        Associative laws &
        Absorption laws \\
        $\mathcal{F} \vee \mathcal{G} = \mathcal{G} \vee \mathcal{F}$ &
        $\mathcal{F} \vee (\mathcal{G} \vee \mathcal{H}) =$
        $(\mathcal{F} \vee \mathcal{G}) \vee \mathcal{H}$ &
        $\mathcal{F} \vee (\mathcal{F} \wedge \mathcal{G}) = \mathcal{F}$ \\
        $\mathcal{F} \wedge \mathcal{G} = \mathcal{G} \wedge \mathcal{F}$ &
        $\mathcal{F} \wedge (\mathcal{G} \wedge \mathcal{H}) =$
        $(\mathcal{F} \wedge \mathcal{G}) \wedge \mathcal{H}$ &
        $\mathcal{F} \wedge (\mathcal{F} \vee \mathcal{G}) = \mathcal{F}$ \\
        \\
        Identity laws & Distributive laws &\\
        $\mathcal{F} \vee \widehat{0} = \mathcal{F}$ &
        $\mathcal{F} \wedge (\mathcal{G} \vee \mathcal{H}) =$
        $(\mathcal{F} \wedge \mathcal{G})  \vee$
        $(\mathcal{F} \wedge \mathcal{H})$ &\\
        $\mathcal{F} \wedge \widehat{1} = \mathcal{F}$ &
        $\mathcal{F} \vee (\mathcal{G} \wedge \mathcal{H}) =$
        $(\mathcal{F} \vee \mathcal{G})  \wedge (\mathcal{F} \vee \mathcal{H})$.
        &
      \end{tabular}
    \end{center}
    The lattice has the additional properties that
    $\mathcal{F} \vee \mathcal{H} = \widehat{1}$ iff
    $\mathcal{F} = \widehat{1}$ or $\mathcal{H} = \widehat{1}$,
    $\mathcal{F} \wedge \mathcal{H} = \widehat{0}$
    iff $\mathcal{F} = \widehat{0}$ or $\mathcal{H} = \widehat{0}$,
    and $b(\widehat{0}) = \widehat{1}$;
  \item
    the map $\mathcal{H} \mapsto b(\mathcal{H})$ is a duality operation:
    $b(b(\mathcal{H})) = \mathcal{H}$;
  \item
    deletion is the dual operation of contraction and vice versa:
    $b(\mathcal{H} \backslash v) = b(\mathcal{H}) / v$ and
    $b(\mathcal{H} / v) = b(\mathcal{H}) \backslash v$;
  \item
    join is the dual operation of meet and vice versa:
    $b(\mathcal{H} \vee \mathcal{F}) = b(\mathcal{H}) \wedge b(\mathcal{F})$ and
    $b(\mathcal{H} \wedge \mathcal{F}) = b(\mathcal{H}) \vee b(\mathcal{F})$;
  \item
    the property of being a minor is preserved by duality:
    $\mathcal{F} \subseteq_{\mathcal{M}} \mathcal{H}$ if and only if
    $b(\mathcal{F}) \subseteq_{\mathcal{M}} b(\mathcal{H})$;
  \item
    minors commute with join and meet:
    \[
    (\mathcal{H} \vee \mathcal{F}) \backslash v =
    \mathcal{H} \backslash v \vee \mathcal{F} \backslash v
    \hspace{20pt}\text{and}\hspace{20pt}
    (\mathcal{H} \vee \mathcal{F}) / v = \mathcal{H} / v \vee \mathcal{F} / v.
    \]
  \end{enumerate}
\end{lemma}

Given a $k$-set $S$, the notation $\mathcal{H} \backslash S$ stands for
$\mathcal{H} \backslash s_1 \ldots \backslash s_k$,
where $\{s_1, \ldots s_k\}$ is an arbitrary ordering of $S$.
The choice of the ordering does not matter by
Lemma~\ref{thm:algebra} (1).
We define $\mathcal{H} / S$ to be $\mathcal{H} / s_1 \ldots / s_k$ likewise.
Finally, given disjoint sets $S$ and $T$
we define $\mathcal{H}[S; T]$ to be $\mathcal{H} \backslash S / T$.
This operation has the nice property
that $\mathcal{H}[S_1 \cup S_2; T_1 \cup T_2] = \mathcal{H}[S_1; T_1][S_2; T_2]$
for any four sets $S_1$, $S_2$, $T_1$ and $T_2$ such that
$(S_1 \cup S_2) \cap (T_1 \cup T_2) = \emptyset$,
which follows from Lemma~\ref{thm:algebra} (1).
We stress the trivial fact that $\mathcal{F} \subseteq_m \mathcal{H}$
iff $\mathcal{F} \cong \mathcal{H}[S; T]$ for some $S$ and $T$.

\section{Outline of the results}
\begin{define}
  [Semi-matching and expanded minor matching]
  Consider a set $\mathcal{S}$ of the form $\{(L_i, S_i)\}_{i=1}^k$,
  a clutter $\mathcal{H}$ and the following conditions:
  \begin{enumerate}
  \item
    $|L_i| = 2$, $L_i \subseteq S_i$ and $S_i \in \mathcal{H}$
    for each $1 \le i \le k$,
  \item
    $\{L_i\}_{i=1}^k$ are pairwise disjoint,
  \item
    [3a.]
    $L_i \not\subseteq S_j$ for $i \neq j$,
  \item
    [3b.]
    $L_i \cap S_j  = \emptyset$ for $i \neq j$,
  \item
    [4.]
    all sets $S \in \mathcal{H}$ contained in $\bigcup_i S_i$
    contain as a subset at least one set $L_i$.
  \end{enumerate}
  We call $\mathcal{S}$ a \emph{semi-matching} if
  it satisfies properties 1, 2, 3a and 4;
  and an \emph{expanded minor matching} if it additionally satisfies 3b.
\end{define}
\noindent
Semi-matchings and expanded minor matchings coincide with induced matchings
in edge clutters.

The following theorem is a generalisation of Alekseev's theorem
from~\cite{alekseev} to clutters.
\begin{thm}
  [Decomposition theorem]
  \label{thm:decomposition}
  For every clutter $\mathcal{H}$,
  $|b(\mathcal{H})|$ is at most the number of semi-matchings of $\mathcal{H}$.
\end{thm}
\noindent
Note the lack of restrictions on the rank of $\mathcal{H}$.

The following lemma establishes that extended minor matchings
and minor matchings are essentially equivalent.
\begin{lemma}
  \label{thm:extended_to_minor}
  If $\mathcal{S} = \{(L_i, S_i)\}_{i=1}^k$
  is an expanded minor matching of a clutter $\mathcal{H}$,
  then $\{L_i\}_{i=1}^k \cong kK_2$ is a minor of $\mathcal{H}$.
  Conversely, if $\{L_i\}_{i=1}^k$ is a $kK_2$ minor of $H$,
  we can find sets $\{S_i\}_{i=1}^k$ such that
  $\{(L_i, S_i)\}_{i=1}^k$ is an extended minor matching of $\mathcal{H}$.
\end{lemma}
\begin{proof}
  To prove the first part,
  contract $(\bigcup_{i=1}^k S_i) \setminus (\bigcup_{i=1}^k L_i)$ and delete
  $V(\mathcal{H}) \setminus (\bigcup_{i=1}^k S_i)$.
  Any non-deleted edge is subsumed by some $L_i$.

  For the second part,
  each $L_i$ must be a contraction of a set $S_i \in \mathcal{H}$.
  It is a routine to verify that
  $\{(L_i, S_i)\}_{i=1}^k$ is an extended minor matching of $\mathcal{H}$.
\end{proof}

Clearly every extended minor matching is a semi-matching.
In the opposite direction, we see that,  in clutters of bounded rank,
semi-matchings contain expanded minor matchings as subsets of linear size.
\begin{thm}
  [Matching theorem]
  \label{thm:matching}
  Suppose $\mathcal{H}$ is a clutter of rank at most $r$
  and $\mathcal{S}$ is a semi-matching of $\mathcal{H}$.
  There is an expanded minor matching
  $\mathcal{S}' \subseteq \mathcal{S}$
  of size at least $|\mathcal{S}|2^{-(r-2)}/(2r-3)$.
\end{thm}
\noindent
Theorem~\ref{thm:main} follows directly from the matching
and decomposition theorems and Lemma~\ref{thm:extended_to_minor}.
Theorem~\ref{thm:matching} is sharp in the sense that if we put $r=2$,
we see that $\mathcal{S} \equiv \mathcal{S}'$.
In this case the entire $\mathcal{S}$ corresponds to an induced matching
(in graph theoretical sense).

The matching theorem does not hold for clutters of unbounded rank.
Consider $\mathcal{H}$ over the vertex set
$\{a_1, \ldots, a_n\} \cup \{b_1, \ldots, b_n\}$
defined by $\mathcal{H} = \{S_i\}_{i=1}^n$,
where $S_i = \{a_i\} \cup \{b_j : 1 \le j \le i\}$.
We see that $\mathcal{H}$ contains a semi-matching
$\mathcal{S} = \{(a_ib_i, S_i)\}_{i=1}^n$ of size $n$,
but $\mathcal{H}$ is $2K_2$-minor-free.

\section{Decomposition and Matching theorems}
We now present the main technical tool to handle semi-matchings.
\begin{define}
  [Expansion]
  Suppose $R_1, \ldots, R_k, C$ are subsets of $\mathbb{N}$
  such that $\{R_i\}$ are pairwise disjoint and $\bigcup_i R_i \subseteq C$.
  We denote the set of functions
  $\func : \{R_i\}_{i = 1}^k \to \bigcup_{i=1}^kR_i$
  such that $\func(R_i) \in R_i$ for each $i$ by $\funcs(\{R_i\}_{i=1}^k)$.
  For each such function $\func \in \funcs(\{R_i\}_{i=1}^k)$ define
  $im_\func := \{\func(R_i) : 1 \le i \le k\}$ to be the image of $\func$.
  For every clutter $\mathcal{H}$ we define
  \[
  \mathcal{H} \circ (\{R_i\}_{i=1}^k, C)
  := \bigvee_{\func \in \funcs(\{R_i\})}
  \mathcal{H}[im_\func; C \setminus im_\func].
  \]
\end{define}
\noindent
Note that $\mathcal{H} \circ (\{R_i\}_{i=1}^k, C) \neq \widehat{1}$
iff $C \setminus im_\func$ is independent
(does not contain any $S \in \mathcal{H}$)
for each $\func \in \funcs(\{R_i\}_{i=1}^k)$,
and hence condition 4 of the definition of semi-matching is equivalent to
$\mathcal{H} \circ (\{L_i\}_{i=1}^k, \bigcup_{i=1}^kS_i) \neq \widehat{1}$.
We normally shorten
$\mathcal{H} \circ (\{L_i\}_{i=1}^k, \bigcup_{i=1}^kS_i)$
to $\mathcal{H} \circ \mathcal{S}$.

\begin{lemma}
  Suppose $\mathcal{H}$ is a clutter, $T \in b(\mathcal{H})$
  and $v \in V(\mathcal{H})$.
  Either $T - v \in b(\mathcal{H} \backslash v)$
  or there is $S \in \mathcal{H}$ and $u \in S$
  such that $T \backslash \{v, u\} \in b(\mathcal{H} \circ (\{\{v, u\}\}, S))$
  and $v \in S - u$.
\end{lemma}
\begin{proof}
  Suppose $T - v \notin b(\mathcal{H} \backslash v) = b(\mathcal{H}) / v$.
  It follows that $v \notin T$ and that there is $T^\prime \in b(\mathcal{H})$
  such that $T^\prime - v \subsetneq T$.
  Let $u \in T \setminus T^\prime$.
  The key observation here is that
  \begin{align}
    \text{ for any } F \in \mathcal{H}
    \text{ either } F \cap T \neq \{u\} \text{ or } v \in F.
    \label{eq:critical_edge}
  \end{align}
  Indeed, if $F \cap (T \cup v) = \{u\}$,
  then $F \cap T^\prime = \emptyset$,
  a contradiction because $T^\prime$ is a transversal.
  Furthermore, we can find a set $S \in \mathcal{H}$
  such that $S \cap T = \{u\}$,
  because $T$ is minimal, and hence $v \in S$ by (\ref{eq:critical_edge}).

  We claim that
  $T - u = T \backslash \{u, v\} \in b(\mathcal{H} \circ (\{\{v, u\}\}, S))$.
  Indeed,
  \begin{align*}
    b(\mathcal{H} \circ (\{\{v, u\}\}, S))
    &= b(\mathcal{H}[u; S-u] \vee \mathcal{H}[v; S-v])\\
    &= b(\mathcal{H})[S-u; u] \wedge b(\mathcal{H}[v; S-v]).
  \end{align*}
  First note that $T - u \in b(\mathcal{H})[S - u; u]$,
  because $T \in b(\mathcal{H})$ and $S \cap T = \{u\}$.
  From (\ref{eq:critical_edge}) we see that $T-u$ intersects all
  sets in $\mathcal{H} \backslash v$,
  hence there is $\widetilde T \in b(\mathcal{H} \backslash v)$
  such that $\widetilde T \subseteq T - u$.
  Furthermore, $\widetilde T \cap S-v = \emptyset$,
  so $\widetilde T \in b(\mathcal{H} \backslash v) \backslash (S-v)$
  $= b(\mathcal{H}[v; S-v])$.
  We deduce that $T - u \in b(\mathcal{H})[S - u; u]$
  and $\widetilde T \in b(\mathcal{H}[v; S-v])$,
  so by the definition of meet ($\wedge$)
  there is a set $R \subseteq (T - u) \cup \widetilde T = T - u$ such that
  $R \in b(\mathcal{H})[S-u; u] \wedge b(\mathcal{H}[v; S-v])$.
  However, if $R \neq T - u$,
  then $R$ is not a transversal of $b(\mathcal{H})[S-u; u]$,
  because $T - u$ is a minimal such transversal,
  in contradiction with the definition of meet.
\end{proof}

\begin{lemma}
  Suppose $\mathcal{H}$ is a clutter,
  $C \in \mathcal{H}$, $R \subseteq C$, $|R| = 2$ and
  $\mathcal{S}' = \{(L_i, S_i')\}_{i=1}^k$ is a semi-matching in
  $\mathcal{H}' = \mathcal{H} \circ (\{R\}, C)$.
  Then there is a semi-matching
  $\mathcal{S} = \{(L_i, S_i)\}_{i=1}^k \cup (\{R\}, C)$ in $\mathcal{H}$,
  such that $\mathcal{H} \circ \mathcal{S} = \mathcal{H}' \circ \mathcal{S}'$
  and $S_i' \subseteq S_i \subseteq S_i' \cup C$ for each $i$.
  Pick one such semi-matching, say the lex-first one,
  and call it $ext(\mathcal{S'}, \mathcal{H}, R, C)$.
\end{lemma}
\begin{proof}
  It follows from the definition of $\mathcal{H} \circ (\{R\}, C)$
  that since $S'_i \in \mathcal{H}'$
  there is at least one set $S_i \in \mathcal{H}$ such that
  $S'_i \subseteq S_i \subseteq S'_i \cup C$ and $R \not\subseteq S_i$
  for each $i \in [k]$.
  Fix an arbitrary such $S_i$ for each $i$ and consider
  $\mathcal{S} = \{(L_i, S_i)\}_{i=1}^k \cup \{(\{R\}, C)\}$.
  From the definitions it immediately follows that $\mathcal{S}$ satisfies
  conditions 1, 2 and 3a. To see that $\mathcal{S}$ satisfies 4 expand
  \begin{align*}
    \mathcal{H} \circ \mathcal{S}
    &= \bigvee_{\func \in \funcs(\{L_i\} \cup \{R\})}
    \mathcal{H}[im_\func; (\cup_i S'_i \cup C) \setminus im_\func]\\
    &= \bigvee_{\func' \in \funcs(\{L_i\})} \bigvee_{r \in R}
    \left(\mathcal{H}[r; C-r][im_{\func'}; \cup_i S'_i \setminus im_{\func'}]
    \right)\\
    &= \bigvee_{\func' \in \funcs(\{L_i\})} \left(\bigvee_{r \in R}
    \mathcal{H}[r; C-r]\right) [im_{\func'}; \cup_i S'_i \setminus im_{\func'}]
    = \mathcal{H}' \circ \mathcal{S}' \neq \widehat{1}.
    &\qedhere
  \end{align*}
\end{proof}

\begin{proof}
  [Proof of the Decomposition Theorem~\ref{thm:decomposition}]
  We denote the set of all semi-matchings of a clutter $\mathcal{H}$
  by $sm(\mathcal{H})$.
  We prove by induction on $|V(\mathcal{H})|$
  that $|b(\mathcal{H})| \le |sm(\mathcal{H})|$.
  If $|V(\mathcal{H})| = 0$,
  then $|b(\mathcal{H})| \le 1$ and $sm(\mathcal{H}) = \{\emptyset\}$.
  Now suppose $v \in V(\mathcal{H})$.
  Let $\mathcal{R}$ be the set of pairs $(R, C)$
  such that $C \in \mathcal{H}$, $R \subseteq C$, $|R| = 2$ and $v \in R$.
  We define $E_{R, C}$ to be the set
  \[
  E_{R, C} :=
  \{ext(\mathcal{S}', \mathcal{H}, R, C) :
  \mathcal{S}' \in sm(\mathcal{H} \circ (\{R\}, C))\}.
  \]
  We see that $E_{R, C} \subseteq sm(\mathcal{H})$
  and $E_{R, C} \cap E_{R', C'} = \emptyset$ for $(R, C) \neq (R', C')$
  as they disagree on the pair of sets containing $v$.
  Furthermore, $sm(\mathcal{H} \backslash v) \subseteq sm(\mathcal{H})$
  and $sm(\mathcal{H} \backslash v)$ is disjoint from each $E_{R, C}$
  because $v \notin L$ for each
  $(L, S) \in \mathcal{S} \in sm(\mathcal{H} \backslash v)$
  and every $\mathcal{S} \in E_{R, C}$ contains a pair, $(R, C)$, including $v$.
  Now we see
  \begin{align*}
    |b(\mathcal{H})|
    &\le |b(\mathcal{H} \backslash v)|
    + \sum_{(R, C) \in \mathcal{R}} |b(\mathcal{H} \circ (\{R\}, C))|\\
    &\le |sm(\mathcal{H} \backslash v)|
    + \sum_{(R, C) \in \mathcal{R}} |sm(\mathcal{H} \circ (\{R\}, C))|\\
    &= |sm(\mathcal{H} \backslash v)|
    + \sum_{(R, C) \in \mathcal{R}} |E_{R, C}| \le |sm(\mathcal{H})|.
    &\qedhere
  \end{align*}
\end{proof}

\begin{proof}
  [Proof of the Matching Theorem~\ref{thm:matching}]
  Suppose $\mathcal{S} = \{(L_i, S_i)\}_{i=1}^\ell$
  is a semi-matching in a rank $r$ clutter $\mathcal{H}$
  and let $G$ be a graph over $[\ell]$
  where $i$ is connected to $j$
  if $|S_i \cap L_j| = 1$ or $|S_j \cap L_i| = 1$.
  It follows that $G$ contains at most $(r-2)\ell$ edges.
  A classic result in graph theory states that
  $\alpha(G) \ge \frac{v(G)^2}{2e(G) + v(G)}$,
  and hence $G$ contains an independent set $I$ of size at least $\ell/(2r-3)$.

  For $\func \in \funcs(\{L_i\}_{i \in [\ell] \setminus I})$ define
  \[
  I'(\func) := \{i \in I :
  \func(L_j) \notin S_i \text{ for each } j \in [\ell] \setminus I\}.
  \]
  Sample $\funcc$ uniformly at random from
  $\funcs(\{L_i\}_{i \in [\ell] \setminus I})$.
  We have $\mathbb{E}|I'| \ge |I|2^{-(r-2)} \ge \ell2^{-(r-2)} / (2\ell-3)$,
  so there must be some $J = I'(\funcc)$ of at least this size,
  where $\funcc \in  \funcs(\{L_i\}_{i \in [\ell] \setminus I})$.

  We claim that $\mathcal{S}_J := \{(L_i, S_i)\}_{i \in J}$
  is an extended minor matching.
  Properties 1 and 2 are inherited from $\mathcal{S}$.
  Property 3b holds because $J \subseteq I$ is stable in $G$.
  Let $\funccc \in \funcs(\{L_i\}_{i \in [\ell] \setminus J})$
  be an arbitrary extension of $\funcc$ and,
  let $C := \bigcup_{i=1}^kS_i$ and $C_J := \bigcup_{i \in J}S_i$.
  We see that $C_J$ is disjoint from $im_{\funccc}$.
  Since every $S \in \mathcal{H}$ with $S \subseteq C$
  contains $L_i$ as a subset for some $i \in [l]$,
  it follows that if $S \subseteq C_J \subseteq C \setminus im_{\funccc}$,
  then $L_i \subseteq S$ for some $i \in J$.
\end{proof}

\section{Applications}
Often instances $I$ of hard problems correspond to a clutter $\mathscr{H}(I)$,
and the solution of $I$ can be read from the blocker of $\mathscr{H}(I)$.
Whether or not $b(\mathcal{H})$ can be computed in polynomial time
from $\mathcal{H}$
with respect to $|\mathcal{H}| + |b(\mathcal{H})|$ is an open problem,
but moreover $|b(\mathscr{H}(I))|$ itself can be exponential.
In this paper we gave a sufficient condition for $|b(\mathcal{H})|$
to be polynomial,
namely $rk(\mathcal{H}) \le r$ and $kK_2 \not\subseteq_m \mathcal{H}$.
Furthermore, $b(\mathcal{H})$ can be computed in polynomial time
with respect to $|\mathcal{H}| + |b(\mathcal{H})|$
for clutters $\mathcal{H}$ of bounded rank \cite{eiter2003new}.

More formally, given a problem $\mathscr{P}$ and a function $\mathscr{H}$
mapping the instances of $\mathscr{P}$ to clutters
and computable in polynomial time, we say that
\emph{the solutions of $(\mathscr{P}, \mathscr{H})$ can be read from the blocker}
if there is an algorithm that given input $I$ and $b(\mathscr{H}(I))$
determines whether or not $I$ is a ``yes'' instance
in time polynomial with respect to the size of its input.

As a first example, suppose $I$ is an instance of the set cover problem,
let $G(I)$ be the hypergraph obtained from $I$ by interchanging
the roles of the vertices and the edges, preserving incidence,
and let $\mathscr H(I) := cl(G(I))$.
We see that the sets of $b(\mathscr H(I))$ correspond to minimal covers,
and hence a minimum sized cover of $I$
can be easily found from $b(\mathscr H(I))$.
A minimum weighted cover can be found the same way,
and more generally, a minimum cover for a monotone oracle $Q$,
that is $Q(S) \le Q(T)$ if $S \subseteq T$,
can be found in polynomial time from $b(\mathscr H(I))$.

Another example is the satisfiability problem (SAT).
On instance $I$ with variables $X$ and clauses $C$
create a clutter $\mathscr H(I)$
with vertices $v_i$ and $\overline{v}_i$ for each variable $x_i \in X$
and a set $S \in \mathscr H(I)$ corresponding to the vertices of $c$
for each clause $c \in C$.
For a satisfying assignment $\sigma : X \to \{0, 1\}$ consider the set
\[
S_\sigma = \{v_i : \sigma(x_i) = 1\} \cup \{\overline{v}_i : \sigma(x_i) = 0\}.
\]
It is clear that $S_\sigma$ is a transversal of $\mathscr H(I)$.
Moreover, every transversal intersecting each pair $\{v_i, \overline{v}_i\}$
at most once can be extended to a transversal of the form $S_\sigma$
for a satisfying assignment $\sigma$.
It follows that $I$ is satisfiable if and only if $b(\mathscr H(I))$
contains a set not containing both $v_i$ and $\overline{v}_i$
for all $x_i \in X$.

The aim is not to give a complete list of applications,
as there is an abundance of problems whose solutions
can be naturally read from the blocker.
For non-trivial applications of transversal enumeration
in artificial intelligence, machine learning, data mining,
model-based diagnosis see \cite{eiter2002hypergraph}.
The solution of each of these problems can be read from the blocker.

Finally, let $\mathcal{C}_{r, k}$
be the class of $kK_2$-minor-free clutters of rank at most $r$.
We conclude that given a problem $(\mathscr{P}, \mathscr{H})$
whose solutions can be read from the blocker and positive integers $r$ and $k$,
there is a polynomial-time algorithm to solve $\mathscr{P}$
for instances restricted to $\mathscr{H}^{-1}(\mathcal{C}_{r, k})$.

The choice for $\mathscr{H}$ is important,
as it affects the size of $\mathscr{H}^{-1}(\mathcal{C}_{r, k})$.
For instance, in the example with satisfiability,
if we had added additional sets $\{v_i, \overline v_i\}$ for each $x_i \in X$,
and asked for a transversal of size $|X|$,
we would have created a large artificial minor matching
in most instances, and hence greatly reduced the size of
$\mathscr{H}^{-1}(\mathcal{C}_{r, k})$.

It is worth noting that for fixed $r$ and $k$
it is possible to test in polynomial time if
$\mathcal{H} \in \mathcal{C}_{r, k}$.
Indeed, suppose $\{(L_i, S_i)\}_{i=1}^k$
is an extended minor matching in $\mathcal{H}$
and $n = V(\mathcal{H})$.
Now let $A = \bigcup_{i=1}^kL_i$, $B = \bigcup_{i=1}^kS_i \setminus A$
and observe that
\[
\mathcal{H}[V(\mathcal{H}) \setminus (A \cup B); B] \cong kK_2;
\]
and that $|A| = 2k$, $|B| \le (r-2)k$.
Therefore, to test if $kK_2 \subseteq_m \mathcal{H}$,
it suffices to test for $O(n^{rk})$ pairs $(S, T)$
if $\mathcal{H}[S; T] \cong kK_2$.
The clutter $\mathcal{H}[S; T]$ can be computed via a naive algorithm
in $O(|\mathcal{H}|^2n)$ time.
\bibliographystyle{alpha}
\bibliography{blocker}
\end{document}